\RequirePackage{fix-cm}
\documentclass[smallextended]{svjour3}       
\smartqed  

\usepackage{graphicx}
\usepackage{lineno,hyperref}
\usepackage{amsmath}
\modulolinenumbers[5]
\usepackage{graphicx}
\graphicspath{{Images1/}}
\usepackage{caption}
\usepackage{subcaption}
\usepackage{epstopdf}
\usepackage[margin=2.5cm]{geometry}
\usepackage{lipsum}
\usepackage{comment}
\usepackage{cite}
\usepackage[round]{natbib}
\usepackage{textcomp}
\renewcommand{\thefigure}{\arabic{figure}}

\begin{document}

\title{The Role of Antigen Competitive Dynamics in Regulating the Immune Response
}
%

\author{
Pantea Pooladvand*  \and Peter S. Kim{$^{\dagger}$} 
\and
Barbara Fazekas de St Groth{$^{\dagger}$} }

\institute{
Pantea Pooladvand \at
              School of Mathematics and Statistics, The University of Sydney, Sydney, NSW 2006, Australia \\
             \email{p.pooladvand@maths.usyd.edu.au}           
            \and
           Peter S. Kim \at
              School of Mathematics and Statistics, The University of Sydney, Sydney, NSW 2006, Australia \\
            \and
              Barbara Fazekas de St Groth \at
            Discipline of Pathology, School of Medical Sciences, Faculty of Medicine and Health, The University of Sydney, Sydney, NSW 2006, Australia
              \and
              * Corresponding author\\
{$\dagger$} Both authors contributed comparably\\
}

\date{Received: date / Accepted: date}

\maketitle

\begin{abstract}
The clonal expansion of T cells during an infection is tightly regulated to ensure an appropriate immune response against invading pathogens. Although experiments have mapped the trajectory from expansion to contraction, the interplay between mechanisms that control this response are not fully understood. Based on experimental data, we propose a model in which the dynamics of CD4+ T cell expansion is controlled through the interactions between T cells and antigen-presenting cells, where T cell stimulation is proportional to antigen availability and antigen availability is regulated through downregulation of antigen by T cells. This antigen-dependent-feedback mechanism operates alongside an intrinsic reduction in cell proliferation rate that may also be responsible for slowing expansion. Our  model can successfully predict T cell recruitment rates into division, expansion and clonal burst size per cell when initial precursors are varied or when T cells are introduced late into an ongoing immune response. Importantly, the findings demonstrate that a feedback mechanism between T cells and antigen presenting cells, along with a reduction in cell proliferation rate, can explain the ability of the immune system to adapt its response to variations in initial conditions or changes that occur later in the response, ensuring a robust yet controlled line of defence against pathogens.
\end{abstract}

\section{Introduction}
\label{intro}

A primary T cell response triggers T cell division and differentiation in a highly organised manner, occurring almost synchronously for a limited time \citep{homann2001differential, beverley2000differences}. Although this controlled expansion is essential to mount an appropriate response to pathogens and protect against autoimmune diseases, the mechanisms that drive this behaviour are not well understood. One explanation is that T cell proliferation and resulting clonal expansion are dependent on continual antigen stimulation, which has been observed experimentally for CD4+ T cells \citep{rabenstein2014differential,yarke2008proliferating}. During an immune response, antigen-specific, naive CD4+ T cells interact with mature antigen-presenting cells (APCs) in lymph nodes and become activated. Here, the theory is that the activated T cells continue to interact with APCs leading to a period of clonal expansion. If T cells are stimulated by antigen to continue dividing, then we must consider the mechanism that limits the expansion. One theory is that T cells downregulate antigen presented by APCs and this behaviour controls the overall T cell clonal expansion, ensuring that larger populations of T cells experience more competition for antigen than smaller populations. Specific antigen downregulation during T cell responses has been observed experimentally \citep{furuta2012encounter, kedl2000t}, and may be mediated via an intercellular molecular transfer process called trogocytosis, reviewed by  \citep{dhainaut2014regulation}.\\

The mechanisms that drive the controlled behaviour of CD4+ T cell proliferation in relation to antigen availability have been explored using mathematical models. Borghans {\em et al.} used ordinary differential equations to experiment with a variety of T cell proliferation functions. Fitting models to {\em in vitro} data providing proliferation estimates of a CD4+ T cell clone over a 24-hour culture, they concluded that T cells compete for antigen sites on APCs, and this competition controls the overall level of T cell expansion \citep{borghans1999competition}.\\

De Boer {\em et al.}, designed a model that assumes CD4+ T cell expansion is proportional to the average concentration of antigen (peptide-MHC concentration) on an APC, and antigen concentration is affected by a T cell downregulation mechanism, called grazing. Grazing (the consumption of surface antigen by interacting T cells) saturates for large T cell numbers under the assumption that APCs can bind to a limited number of T cells at a time \citep{de2013antigen}. This model successfully explained {\em in vivo} experimental results by Quiel {\em et al.}, in which initial T cell numbers were varied and the factor of expansion (the ratio of the number of T cells at day 7 to the number of T cells at day 1) was measured \citep{quiel2011antigen}. Following this work, Mayer {\em et al.} also explained the experimental results by Quiel {\em et al.} using a T cell competition model in which T cell proliferation is regulated by antigen presentation (concentration of peptide-MHC molecules), the number of T cells present and the affinity of T cells to antigen. In this model, T cells do not influence the downregulation of peptides, and overall expansion is limited by antigen decay \citep{mayer2019regulation}.\\

The role of CD4+ T cells has also been included in models capturing the broader dynamics and interactions involved in cell-mediated immunity response. Pappalardo {\em et al.} adopted an agent-based approach to describe the response of immune cells (B lymphocytes, CD4+ T cells, CD8+ T cells, macrophages and dendritic cells) to influenza A virosome when administered in conjunction with pre-selected adjuvants (a substance that can improve the immune response) \citep{pappalardo2016computational}. In this model CD4+ T cells activate and proliferate through interactions with B cells and macrophages. The model predicts the most suitable adjuvant to be used in conjunction with influenza A virosome based on the response of antibodies (immunoglobulin G). The authors show that the model prediction is in good agreement with in vivo results. Pennisi {\em et al.} modelled the efficacy of the RUTI\textsuperscript{\sffamily\textregistered} vaccine against Mycobacterium tuberculosis (MTB) by using the Universal Immune System Simulator (UISS), an agent-based computational framework \citep{pennisi2019predicting}. The authors used this simulator to build a complex system of interactions, considering both innate and adaptive immunity. The model captures interactions between cells and cytokines released during the immune response. Naïve CD4+ T cells differentiate into phenotypes and interact with B cells driving the response. The model is in good agreement with the results from a phase II clinical trial where the subjects with latent tuberculosis infection were treated with RUTI\textsuperscript{\sffamily\textregistered}. 

\begin{figure}
    \centering
    \includegraphics[width=5cm]{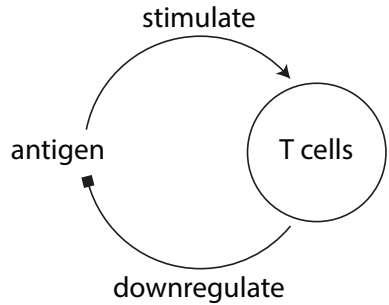}
    \caption{T cell stimulation and regulation model. Presence of antigen stimulates T cells into activation and division. Naive and activated T cells downregulate antigen during this process.}
    \label{fig:flow}
\end{figure}

In this paper, we consider a model in which CD4+ T cell expansion is controlled at a rate proportional to antigen concentration and antigen concentration is downregulated by T cells as in Figure~\ref{fig:flow}; however, our model does not assume saturation due to limitations in space or binding of T cells to MHC molecules. We base this assumption on experiments by Spencer {\em et al.} \citep{spencer2020antigen}, where the behaviour of monoclonal, high affinity CD4+ T cells responding to specific peptide-MHC complexes {\em in vivo} is observed. The authors showed that the response to a second, independent antigen did not limit the initial response even when the same dendritic cells presented both antigens (refer to Figure 4-7 in~\citep{spencer2020antigen}). Our model also considers an intrinsic reduction in cell proliferation upon subsequent divisions, which is implied by fitting the model to experimental data. Building a mechanistic system of delay differential equations, we track T cell divisions and total T cell numbers to explain the dynamics of T cell expansion during a primary immune response from two separate experimental scenarios: one in which the overall clonal expansion of CD4+ T cells is measured against initial T cell numbers, and another in which T cell recruitment and division is measured for a cohort of T cells that is introduced later in an ongoing immune response. By fitting our model to these scenarios we demonstrate that the regulation of T cell clonal expansion may result from an antigen competitive environment along with a reduction in T cell proliferation rate. 

The work is organised as follows. In Section~\ref{Exp overview} we summarise the experimental method. In Sections~\ref{Model} and \ref{Param est} we introduce our model and discuss model calibration to three experimental data sets. In Section~\ref{ResultsSection} we compare our model prediction to the experimental results. Finally, we discuss the implications of our model and future work in Section~\ref{Discussion}.

\subsection{Overview of experimental method}\label{Exp overview}
We summarise the experiments of Smith {\em et al.} and Spencer {\em et al.}, which we will use to quantify the dynamics of CD4+ T cell expansion during an immune response. 
\begin{itemize}
    \item Experiment 1 \citep{Smith2020.09.09.290627} measures T cell clonal expansion during an immune response in relation to initial T cell numbers. In this experiment, transgenic lymph node T cells are injected into non-transgenic Ly5.1 congenic recipient mice in doses of $10^4$, $10^5$, $10^6$ and $10^7$ cells. T cells in the system are measured on day 0 and all mice receive a fixed dose (10$\mu$g) of cognate antigen, injected on the same day. T cells are counted on days 7 and 42. Results in Figure~\ref{fig:Experiment1} show the number of Tg+ T cells per $10^5$ leukocytes measured in the draining lymph nodes.\\
    \item Experiment 2 \citep{spencer2020antigen} measures the percentage of T cells recruited and their division profiles when CD4+ T cells are injected later in an ongoing immune response. In this experiment, two cohorts of 5C.C7 T cells are injected into recipient mice. The first cohort is the competing T cells and the second is the CFSE-labelled cells, which we are tracking. Group (i) is the control group where no competing cohort is injected. In Group (ii), both cohorts are injected on the same day, and in Group (iii), the competing cohort is injected 24 hours prior to the labelled cells. All recipient mice receive the same dose of cognate antigen ($10\mu$g) 12 hours prior to the start of experiment. The results for this experiment are shown in Figure~\ref{fig:Experiment2}.\\
    \item Experiment 3 \citep{spencer2020antigen} is similar to Experiment 2 but with an extended delay for the CFSE-labelled cells. In this experiment, we have a control group Group (i) with no competing cohort, Group (ii) with a 24-hour delay between the competing cohort and the labelled cells and Group (iii) with a 72-hour delay for the labelled cells. Results are shown in Figure~\ref{fig:Experiment3}.
\end{itemize}
We optimise the parameters in our model by fitting to all three experiments simultaneously.

\section{Model}\label{Model}
We consider a scenario in which adoptively transferred, naive CD4+ T cells interact with mature antigen-presenting cells (APCs) in lymph nodes during a primary immune response. For simplicity, we model the general level of antigen, instead of explicitly introducing a population of APCs.  Our model includes antigen, $A(t)$; naive T cells, $N(t)$; and activated T cells, $T_1(t)$, $T_2(t)$, \dots, $T_i(t)$, \dots, where the subscript $i$ represents the number of cell divisions undergone by the population by time $t$. Along with estimating the total number of T cells in time, we will also use this model to estimate the number of divisions in the timespan of 3.5 and 5.5 days. These are the timespans for Experiment 2 and 3 respectively. For this reason we use a system of delay differential equations to capture the time delay between divisions more accurately. The following system of delay differential equations (DDE) describes the interactions between our populations:
\begin{align}
 A'(t) = & f_A(t) - \left(s_N N + \sum_{i=1} ^{\infty} s T_i \right)A - d_A A, \label{eq:A}\\
 N'(t) = & f_N(t) - r_N N A, \label{eq:N} 
\end{align}

\noindent with activated T cell populations
\begin{align}
 T_1'(t) & = 2r_N N(t-\sigma) A(t-\sigma) - r_1T_1 A - dT_1, \label{eq:T1}\\
 T_2'(t) & = 2r_1e^{-d\tau}T_1(t-\tau) A(t-\tau) - r_2T_2 A - dT_2, \label{T2}\\
 & \vdots \nonumber \\
 T_i'(t) & = 2r_{i-1} e^{-d\tau}T_{i-1}(t-\tau) A(t-\tau) - r_{i}T_i A - d T_i, \label{eq:Ti}\\
 \vdots  \nonumber 
\end{align}

 In equation (\ref{eq:A}), antigen is supplied to the lymph nodes at rate $f_A(t)$. We assume  downregulation of antigen occurs at mass-action rates proportional to antigen availability and T cell population size. We also consider that downregulation may occur at different rates for naive cells compared with activated cells. Therefore the rates $s_N NA$ and $s T_i A$ describe the downregulation of antigen by naive and activated T cells respectively. Antigen decays at rate $d_AA$.\\

In equation (\ref{eq:N}), naive T cells, $N(t)$, enter the lymph nodes at rate $f_N(t)$ and activate at mass-action rate $r_N NA$ upon interaction with antigen with coefficient $r_N$. The half life of naive CD4+ T cells in mice is approximately 7 weeks \citep{den2012maintenance}, so for the timescale of our model, we assume naive T cells do not decay and remain constant in lymph nodes if no antigen is present.\\

\begin{figure}
    \centering
    \includegraphics[scale=0.35]{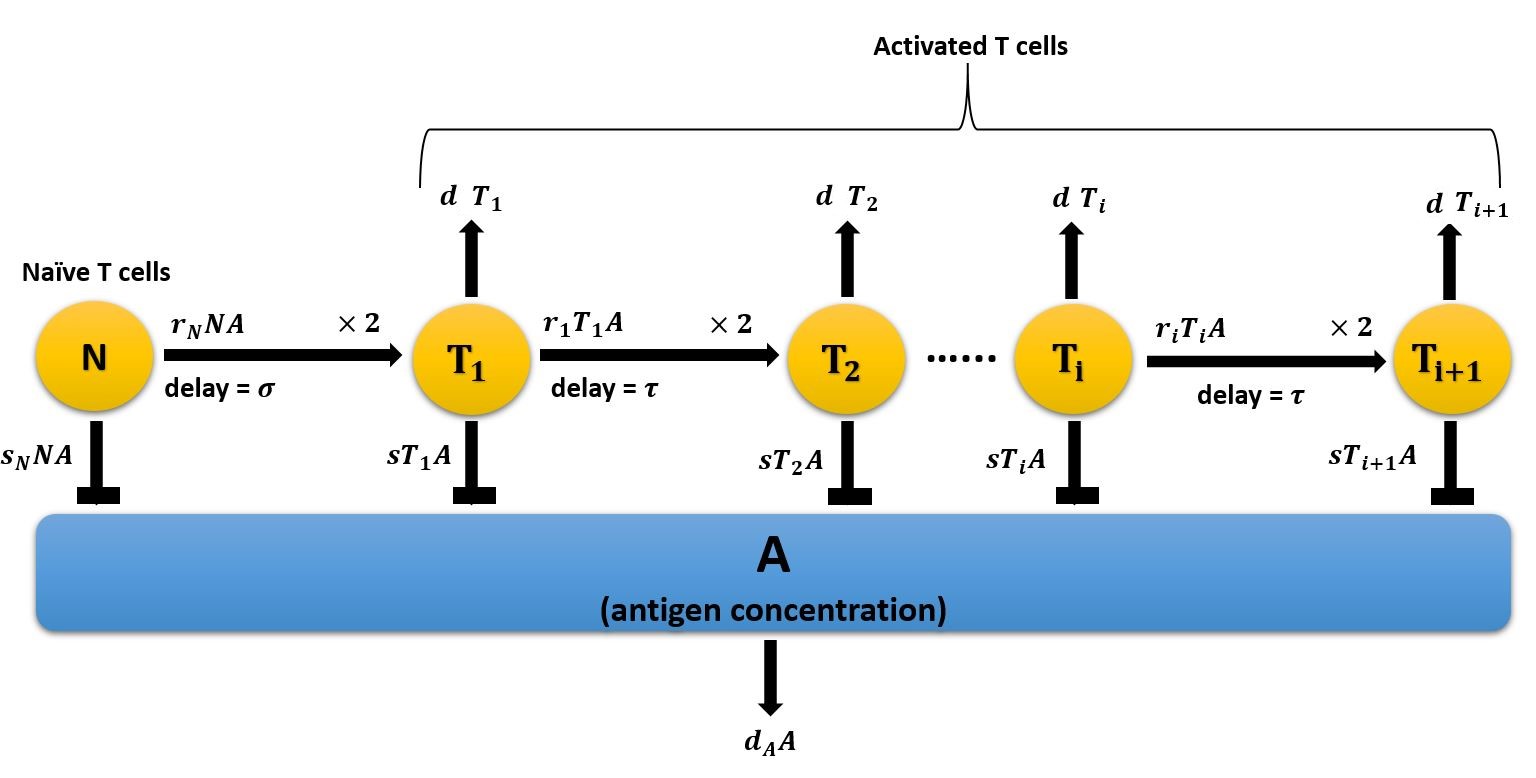}
    \caption{T cell dynamics as described by the model. We implicitly model the interactions of T cells with APCs by measuring a general concentration of antigen, $A(t)$, available in lymph nodes. Naive T cells, $N(t)$, downregulate antigen at rate $s_N NA$ and activate at rate $r_N NA$, proliferating $\sigma$ hours later. The daughter cells move to state $T_1$. Cells in state $T_1$ decay at rate $d T_1 $ and proliferate at rate $r_1T_1A$, moving to state $T_2$ a time delay of $\tau$ hours later. Cells continue to proliferate, decay and downregulate antigen in this manner. A cell in state $T_i$ proliferates at rate $r_i T_i A$, decays or exits lymph nodes at rate $d T_i$ and downregulates antigen at rate $s T_i A$. In addition, antigen decays at rate $d_A A$.}
    \label{fig:flow3}
\end{figure}

The first term of equation (\ref{eq:T1}) accounts for T cells that have undergone their first round of divisions arriving in state $T_1$. The time delay $\sigma$ accounts for the time taken for naive T cells to activate and divide once. The next round of divisions occur at rate $r_1 T_1 A$ and daughter cells move to state $T_2$ after a delay of $\tau$. The population of T cells continue  to divide in this manner with each T cell state having a distinct activation/proliferation rate $r_i$. Equation (\ref{eq:Ti}) describes the general form for the number of T cells in state $i$. We note that subsequent divisions following $T_1$ experience the same delay of length $\tau$. The factor $e^{-d \tau}$ accounts for the cells that decay or leave the lymph nodes during the delay period. All T cells $T_i$ decay or exit lymph nodes at rate $dT_i$. The population interactions are illustrated in Figure~\ref{fig:flow3}.\\

\section{Parameter estimates}\label{Param est}

In keeping with Experiment 1's design, we scale the T cell population as numbers of cells per $10^5$ leukocytes and antigen as proportions of the injected dose $A_{inj}$, so that the initial dose scales to 1. We administer this dose using the supply rate $f_A(t)$. As we have a system of delay differential equations, we need to specify a history for our populations, so we set all populations, $A(t)$, $N(t)$ and $T_i(t)$ to zero, for $t<0$. Initial populations at time 0 will be outlined in Section~\ref{ResultsSection} as they are unique to each experiment.\\

\begin{figure}
    \centering
    \includegraphics[width=\linewidth,trim={2cm 8cm 2cm 5cm},clip]{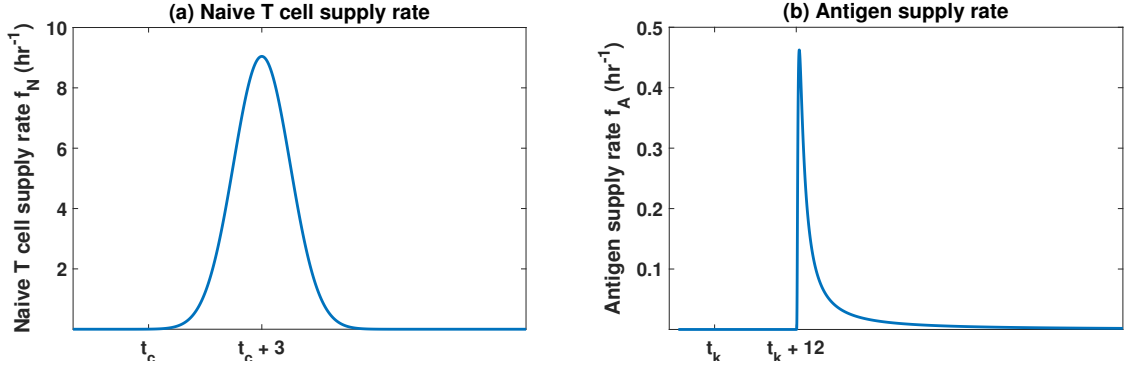}
    \caption{\textbf{(a)} Naive T cell supply function, $f_N(t)$. We use a normal distribution function with cells injected at time $t_c$ and peak numbers arriving in lymph nodes 3 hours following injection time. \textbf{(b)} Antigen supply function, $f_A(t)$, using a long-tailed distribution. Antigen begins to enter the lymph nodes in large quantities 12 hours following injection time $t_k$. The rate declines but small levels of antigen continue to flow into the lymph nodes.}
    \label{fig:Tcell_antigen_supply}
\end{figure}

After transfer, CD4$^+$ T cells migrate to lymph nodes within minutes with substantial numbers entering 2 to 4 hours later \citep{liou2012intravital}. Using this information, we set a normal distribution function for the naive T cell supply at time $t$:
\begin{equation}
f_N(t) \sim \mathcal{N}(\mu,\,p^{2}),
\label{eq:fN}
\end{equation}
where $\mu = t_c + 3$ , $t_c$ is the time of injection of naive T cells and $p = 3/4$ as illustrated in Figure~\ref{fig:Tcell_antigen_supply}a. Our choice in $\mu$ and $p$ means that the supply peak is at 3 hours following T cell transfer, with the majority (91\%) of cells supplied by the fourth hour. Increasing $p$ up to 5 so that T cells are supplied over 48 hours does not significantly alter the outcome.\\

Previous experiments in mice have shown that activated dendritic cells begin to appear in the draining lymph nodes 8 to 16 hours following footpad injections \citep{liou2012intravital}. Allan {\em et al.} measured small numbers of DCs in the lymph nodes at 12 hours with peak numbers 1 to 2 days after treatment~\citep{allan2006migratory}. Data presented by Tomura {\em et al.} described a steep rise in DCs for the first 24 hours following treatment with a very small increase in the next 24 hours~\citep{tomura2014tracking}. With this information in mind, we assume that antigen enters the lymph nodes 12 hours following injection time, $t_k$, with the majority of antigen supplied within the first day. We use a stable distribution
\begin{equation}
f_A(t) \sim \mathcal{S}(\alpha,\beta,\gamma, \delta;0)
\label{eq:fA}
\end{equation}
with shape parameters $\alpha = 0.5$ and $\beta = 1$, scale parameter $\gamma = 1$ and location parameter $\delta = t_k+12$, to supply the lymph nodes with antigen. Antigen injection time, $t_k$, depends on the experiment and will be detailed in Section~\ref{ResultsSection}. These parameters produce a long tailed supply function as in Figure~\ref{fig:Tcell_antigen_supply}b. Using a decreasing step function is also suitable for this system and produces similar results as long as we are in the range of 2 days for the delivery of the majority of antigen.\\

The first division of CD4+ T cells takes approximately 24 hours \citep{obst2015timing}, so we let $\sigma = 24$ hours. In recreating Experiment 2 and Experiment 3, we need to track the number of divisions experienced by T cells in a short time period of 3.5 and 5.5 days respectively. This means that the initial delay is an important parameter to capture the wait time between T cell activation and the time at which T cells begin to rapidly divide and downregulate antigen. The equivalent ordinary differential equation system without delay cannot capture the initial mean division time even at very low proliferation rates. The delay of subsequent divisions $\tau$ is fitted to experimental data with an optimum fit value of 3.9863 which agrees with the observed range of 4 to 6 hours \citep{rabenstein2014differential, van2016proliferating}.\\

\begin{figure}
     \centering
         \includegraphics[scale = 0.3]{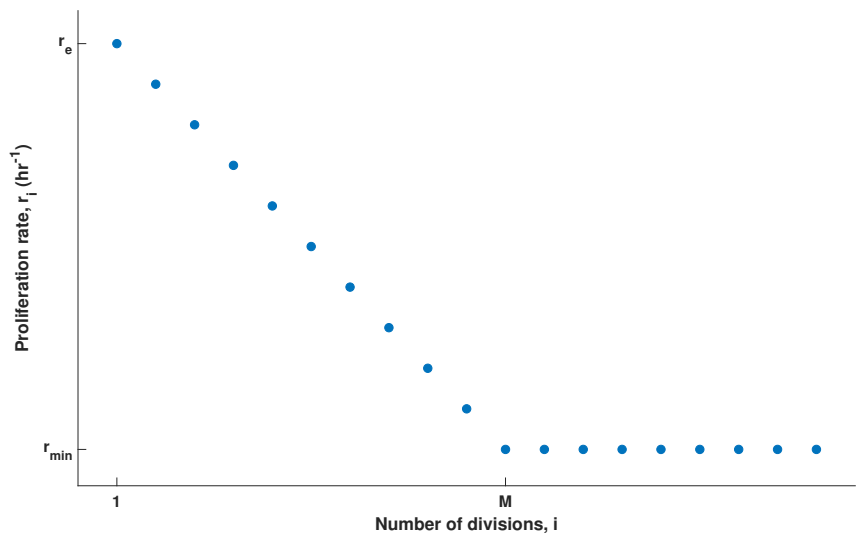}
         \caption{T cell proliferation rate, $r_i$. We assume T cell proliferation decreases in a linear manner with each division until $M$ divisions. Following the $M$th division, $r_i$ reaches a constant rate $r_{min}$.}
        \label{fig:ri and di}
\end{figure}

We model T cell activation and division as proportional to antigen availability. As stated previously, we do not consider saturation or space limitations as a cause for reduction in T cell divisions; instead, we assume that T cell division rate slows with subsequent divisions. This reduction in proliferation rates has been observed experimentally. For example, Homman {\em et al.} measured the division rates of CD4+ T cells during an immune response and concluded that the average division rates are 9-20 hours between days 5 and 7 reducing to 20-40 hours between days 7 and 9~\citep{homann2001differential}. We will treat this reduction in division as a function of the number of divisions $i$, rather than a reduction in time. We choose a proliferation function for $r_i$ that is linearly decreasing with each division, reaching a constant division rate $r_{min}$ after $M$ divisions as in Figure~\ref{fig:ri and di}. Therefore our proliferation parameter $r_i$ is a piece-wise function of the number of divisions $i$:

\begin{equation}
    r_i = 
    \begin{cases}
    r_e(1-g(i-1)), & \text{if } 1 \leq i < M,\\
    r_e(1-gM),     & \text{if } i \geq M,
    \end{cases}
    \label{eq:ri}
\end{equation}
where $r_e$ is the maximum division rate of T cells and $r_{min} = r_e(1-gM)$. We note that an exponentially decaying proliferation function such as $r_i = r_e(0.85)^i$ or a step function constant at $r_i = r_e$ when divisions are below approximately 7 or 8 and reducing to a fraction of $r_e$ for larger divisions describes the data equally well. Therefore, the model does not heavily depend on the type of function, rather that the proliferation rate decreases as cells divide. The experimental results in Figures~\ref{fig:Experiment2}c and~\ref{fig:Experiment3}c show that T cells divided up to 7 times in 2.5 days following an immune response, so we choose $M = 10$. Our choice in using a linear function is motivated so as to not over-dampen the divisions experienced per cell. We fit $r_e$ and $g$ to the experimental data. The fitted maximum proliferation rate $r_e = 1.5412$. This value equates to an average doubling time of $\frac{1}{1.5412} + 4 \approx 4.6$ hours in agreement with experimental results reporting 5 to 8 division in 24 to 36 hours \citep{van2016proliferating} and another reporting division every 4 to 6 hours \citep{rabenstein2014differential}.

We also need to consider the activation and first division rate, $r_N$, of naive T cells. This parameter controls the number of cells that are recruited to the response and is measured in Experiment 2 (see Figure~\ref{fig:Experiment2}b) and Experiment 3 (see Figure~\ref{fig:Experiment3}b). We fit this parameter to the experimental data.

Given that we are implicitly modelling presentation of antigen on APCs, we estimate that the rate $d_A$ corresponds to the turnover of APCs. Mature dendritic cells have an estimated half life of 2 to 3 days under steady-state conditions~\citep{diao2006situ, kamath2002developmental}. Assuming the majority of our antigen-presenting cells are dendritic cells, we use 2.5 days (60 hours) for the half life of our $A(t)$ population so that $d_A = \ln (2)/60 \approx 0.01 \text{hr}^{-1}$.\\

Lambrecht {\em et al.}, showed that CD4+ T cells begin to exit the lymph nodes approximately 2 days following immunisation \citep{lambrecht2000induction}. From our simulations, cells enter $T_1$ state approximately 1.5 days following infection. For simplicity, we fit a constant clearance rate $d$ to all populations of T cells. The clearance rate accounts for T cells leaving the lymph nodes in the early stages of the response and dying at later stages. For the purpose of this model, we consider these cells as permanently removed from the lymph nodes, although it is likely that a small number will continue to recirculate through the priming node.\\

Lastly, we consider downregulation of antigen for naive and activated T cells with rates $s_N$ and $s$, respectively. By optimising $S_N$ to data, we find that the optimal value for this parameter is very small and the confidence interval contains the null value. Therefore, we simplify our model by setting this parameter to zero. The removal of this parameter means that we assume naive T cells do not significantly contribute to the downregulation of antigen. The parameter $s$ is fitted to experimental data. \\

All parameter estimates can be found in Table~\ref{Tab:Parameter estimates}. Parameters $r_e$, $r_N$, $g$, $d$, $\tau$ and $s$ are optimised to data by solving model (\ref{eq:A}-\ref{eq:Ti}) in MATLAB R2018a and performing a simultaneous fit of model output to experimental data in Figure~\ref{fig:Experiment1}, Figure~\ref{fig:Experiment2}b, Figure~\ref{fig:Experiment2}c, Figure~\ref{fig:Experiment3}b and Figure~\ref{fig:Experiment3}c. We use the built-in fitter `lsqcurvefit', nonlinear least-squares, trust-region-reflective algorithm.\\

\begin{table}[h]
\renewcommand{\arraystretch}{1.5}
\renewcommand*\thetable{\Roman{table}}
    \centering
\begin{tabular}{|p{1.5cm}|p{5cm}|p{1cm}|p{3cm}|p{2cm}|}
\hline
\multicolumn{5}{|c|}{Parameter estimates} \\
\hline
 Parameter & Description & Units & Estimate & Source/Conf Interval\\
 \hline
 $r_e$  & Maximum T cell proliferation rate & hr$^{-1}$ & 1.5412  & $1.5333-1.5491$\\
 \hline
 $r_N$  & Naive T cell proliferation rate & hr$^{-1}$ & 0.0497 & $0.0443-0.0551$ \\
 \hline
 $g$  & Change in T cell proliferation rate per division & - & 0.0994 & $0.0899-0.1089$ \\
 \hline
 $d$ & T cell clearance rate & hr$^{-1}$  & 0.0009 & $0.0003-0.0016$  \\
 \hline
 $M$ & Largest number of divisions before $r_i$ reaches a constant  & - & 10 & Estimated\\
 \hline
$s$ & Antigen downregulation & hr$^{-1}$ & 0.0009 & $0.0007-0.0011$\\
\hline
$\sigma$ & Naive T cell activation/first division delay & hours & 24 & \citep{obst2015timing} \\
 \hline
$\tau$ & Activated T cell proliferation delay & hours & 3.9796 & $3.9752-3.9841$  \\
 \hline
$d_A$ & Death rate of APCs & hr$^{-1}$ & 0.01 & \citep{diao2006situ, kamath2002developmental} \\
\hline
$t_c$ & Time of injection for naive T cells & - & Details in Section~\ref{ResultsSection} & \\
\hline
$t_k$ & Time of injection for antigen & - & Details in Section~\ref{ResultsSection} & \\
\hline
$f_A(t)$ & Supply rate of antigen & hr$^{-1}$ & $\mathcal{S}(0.5,1,1,t_k+12;0)$ & \citep{allan2006migratory, tomura2014tracking}\\
\hline
$f_N(t)$ & Supply rate of naive T cells & hr$^{-1}$ & $\mathcal{N}(t_c+3,\frac{3}{4}^{2}$) & \citep{liou2012intravital}\\
\hline
\end{tabular}
 \caption{List of parameter estimates. Parameters $r_e$, $r_N$, $g$, $d$, $\tau$ and $s$ are optimised by solving model (\ref{eq:A}-\ref{eq:Ti}) numerically in MATLAB, using 'dde23' and fitting the model output simultaneously to Experiments 1, 2 and 3 (see Figures~\ref{fig:Experiment1}, \ref{fig:Experiment2} and \ref{fig:Experiment3}) using the built-in nonlinear least-squares, trust-region-reflective algorithm `lsqcurvefit'.}
  \label{Tab:Parameter estimates}
\end{table}


\section{Results}\label{ResultsSection}
We solve (\ref{eq:A})-(\ref{eq:Ti}) numerically using `dde23' in MATLAB R2018a. We first consider the case in which initial naive T cell counts are varied for a fixed dose of antigen, see Figure~\ref{fig:Experiment1}. In this experiment, naive CD4+ T cells are injected into recipient mice in four different doses, increasing incrementally by factors of 10. The T cells are injected 3 days prior to the start of the experiment and numbers are measured at day 0. We run the model with each initial condition $N_0$ as the average measured at day 0, namely $N_0 = 0.1$, $1.3$, $8.5$ and $94.7$ (cells/$10^5$ leukocytes). In the experiment, no additional T cells are introduced into the system, so $f_N(t) = 0$ for all $t$. Since the cognate antigen is injected on day 0 of the experiment, we set $t_k = 0$ for the antigen supply function, $f_A(t)$ described in equation (\ref{eq:fA}). To track total T cells in time, we calculate the sum of T cells in all states at time t as
\begin{equation}
 \text{total T cells} = N(t) + \sum^{20}_{i=1} T_i(t) + D_N(t) + \sum^{19}_{i=1} D_i(t). \label{eq:Total T cells}  
\end{equation} 
In this equation, $D_N$ represents cells undergoing first division, where
\begin{equation}
    D_N(t) = \int \left(r_N NA - r_NT_N(t-\sigma)A(t-\sigma) \label{eq:D_N} \right)\ dt, 
\end{equation}
and $D_i$ represents cells undergoing the $i^{th}$ division, where
\begin{equation}
    D_i(t) = \int \left(r_i T_i A - r_i T_i(t-\tau)A(t-\tau)e^{-d\tau} -d D_i \right)\ dt,  \label{eq;D_i}
\end{equation}
so that $D_N$ is the transition state between $T_N$ and $T_1$, and $D_i$ is the transition state between $T_i$ and $T_{i+1}$. From equation (\ref{eq:D_N}), naive T cells enter at rate $r_N NA$ and remain in this state for $\sigma$ hours before exiting to $T_1$. For higher division states $D_i$, the exponential factor in the second term accounts for the fraction of cells that survive the division period of $\tau$ hours, consistent with equations (\ref{eq:T1})-(\ref{eq:Ti}). We truncate simulations at 20 divisions as number of dividing cells at this threshold drop below $10^{-5}$. \\

\begin{figure}
    \centering
    \includegraphics[width=\textwidth]{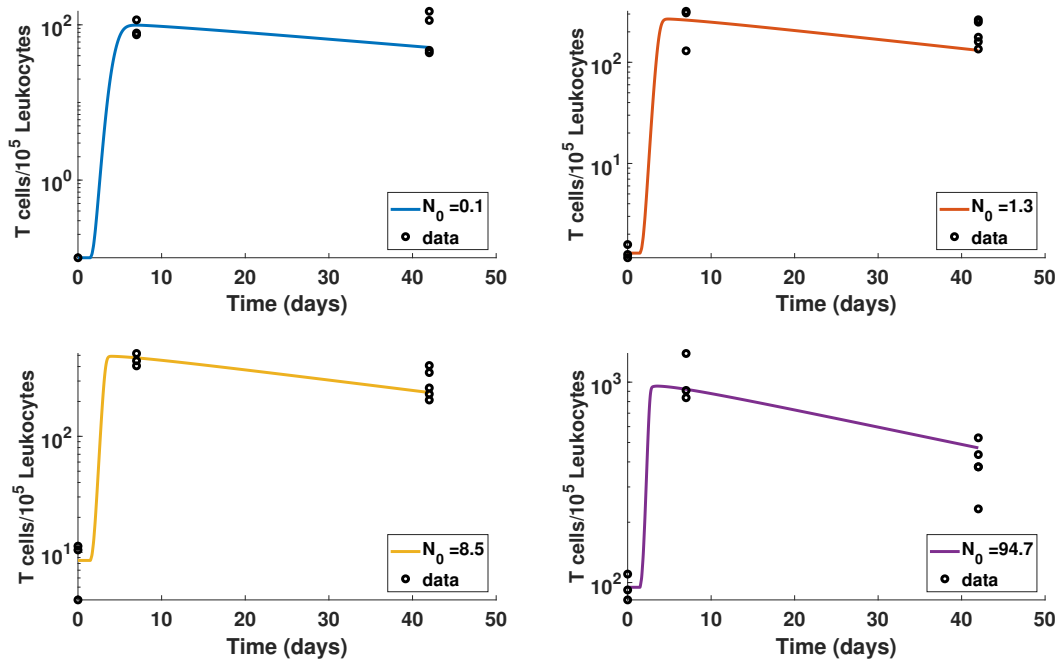}
    \caption{Model fit to data. Total T cells in time presented on a semi-log scale. T cells injected 3 days prior to start of experiment and measured at day 0. Data points represent the number of T cells measured on days 0, 7 and 42. Lines represent model prediction for each precursor dose over 42 days.}
    \label{fig:fit_to_data1}
\end{figure}

\begin{figure}
    \centering
    \includegraphics[width = \textwidth, trim={2cm 2cm 2cm 2cm},clip]{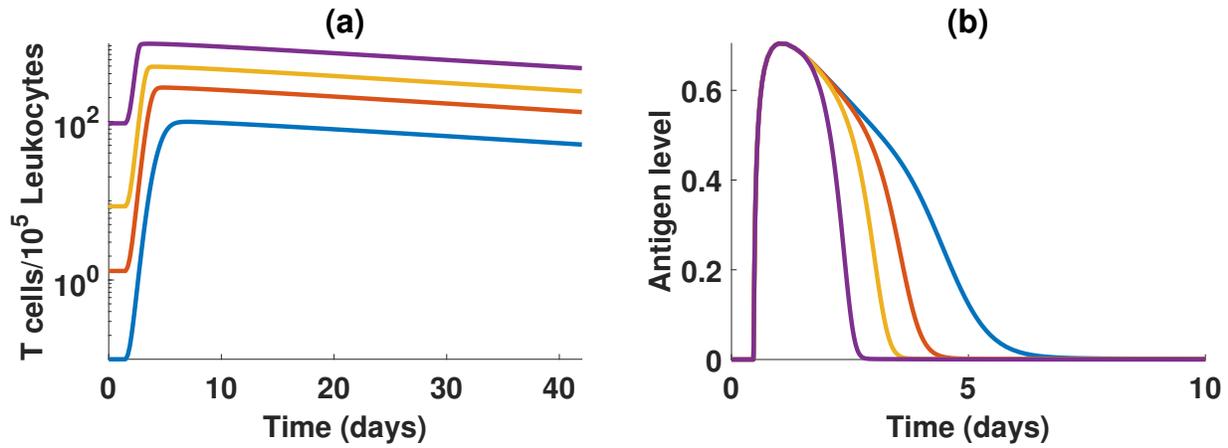}
    \caption{Comparisons between precursor doses of T cells. \textbf{(a)} Tracking total T cells for the 3 measured time points from data shows that T cell expansion decreases with increasing precursors. Fold difference, between smallest and largest precursors of 1000 at day 0 reduces to 9.7 on day 7 followed by 5.0 on day 42. \textbf{(b)} Antigen availability in time reduces as precursors are increased due to downregulation of antigen by T cells.}
    \label{fig:divisions1}
\end{figure}

In Figure~\ref{fig:fit_to_data1}, we have the total number of T cells over 42 days from expansion to contraction (semi-log) fitted to experimental data taken at days 0, 7 and 42. The results show that our model can capture the behaviour in the data for varying initial precursor doses of antigen-specific naive T cells. If we compare the rise in T cell numbers between the precursors at the three time points, 0, 7 and 42, as in Figure~\ref{fig:divisions1}a, there is a clear reduction in clonal expansion of T cells as the precursor dose is increased from $N_0 = 0.1$ (blue) to $N_0 = 94.5$ (purple) in agreement with results in Figure~\ref{fig:Experiment1}.
As initial precursors are increased, antigen levels decrease more rapidly due to T cell downregulation of antigen, controlling the size of clonal expansion. In Figure~\ref{fig:divisions1}b, we see that antigen levels diminish faster as initial precursors rise. In the smaller precursors, where antigen is available over a longer period of time, the cell intrinsic proliferation rate, $r_i$ limits the expansion. Together, antigen availability and reduction in cell proliferation, capture the observations from Experiment 1.

\begin{figure}
    \centering
    \includegraphics[scale = 0.3]{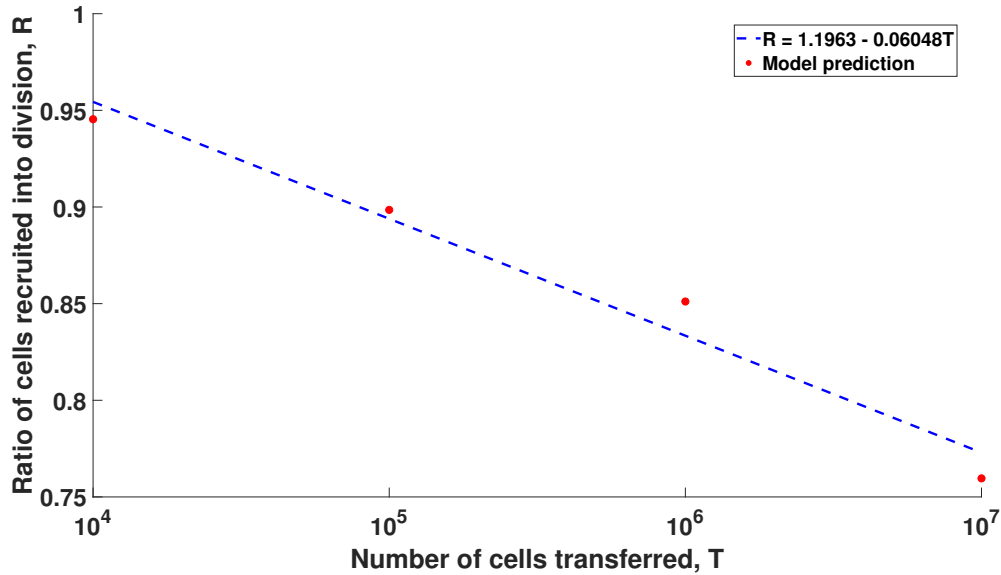}
    \caption{Predicted ratio of recruited T cells to transferred T cells. Model estimates a $6\%$ decrease in recruitment of T cells to division for each factor of 10 increase in initial T cell counts. R-Squared = 0.969.} 
    \label{fig:regression_cells_to_activation}
\end{figure}

The fold difference (ratio of largest population ($N_0 = 94.7$) to smallest population ($N_0 = 0.1$)) from experimental results are 1000 at day 0 and 12 at day 7. Our model predicts a fold difference of 9.3 on day 7. We can also use our model to estimate the percentage of cells that are recruited into division by calculating the ratio of naive cells on the last day of experiment to naive cells on day 0 as
\begin{equation}
    \textbf{\% of cells recruited into division} = \left(1-\frac{N(\text{last day of experiment})}{N(0)} \right) \times 100. \label{eq:recruitment}
\end{equation} 
In Figure~\ref{fig:regression_cells_to_activation} we have the ratio of recruitment into division for each group of transferred cells. From the line of best fit we can estimate a reduction of $6\%$ per factor 10 increase in initial naive cell counts.\\

In the next experiment, the efficiency of CD4+ T cells to contribute to an ongoing immune response is determined by transferring two separate cohorts of naive T cells into recipient mice and increasing the time delay between the two transfers, Figure~\ref{fig:Experiment2}a. The results in Figure~\ref{fig:Experiment2} track the delayed population (CFSE-labelled cells) and compare the effect of delay on the percentage of cells recruited into division, Figure~\ref{fig:Experiment2}b, and division profiles on the final day of experiment (day 3.5), Figure~\ref{fig:Experiment2}c. To simulate these results, we introduce a second cohort of T cells with naive cells, $\Tilde{N}(t)$, and activated cells in state $i$, $\Tilde{T}_i$. We also define a T cell supply rate $f_{\tilde{N}}(t)$, such that $f_{\tilde{N}}(t) \sim \mathcal{N}(\tilde{\mu},p^2)$, where $\tilde{\mu} = \tilde{t}_c+3$, $p = 3/4$ and $\tilde{t}_c$ is the transfer time for the second cohort of T cells. These populations mirror our original system, (\ref{eq:A}-\ref{eq:Ti}), with total T cells as measured in equation (\ref{eq:Total T cells}).
We will refer to our original population as the CFSE-labelled cells (or just labelled cells) and the second cohort as the competing cells.
In this experiment, antigen is injected approximately 12 hours prior to the first injection of naive T cells, so we set time of antigen injection $t_k=-12$ for our antigen supply function $f_A(t)$. Since the injected dose of antigen in this series of experiments is the same as the previous, we assume that a similar antigen concentration will be presented by the APCs, and therefore, total antigen dose is 1. 
From experimental design in Figure~\ref{fig:Experiment2}a, in all three groups (i) to (iii), our labelled T cell population $N(t)$ is injected 24 hours after day 0. Therefore we set time of injection $t_c = 24$ and $\mu = 27$ for our T cell supply function $f_N(t)$. Each injected dose of T cells is $2 \times 10^6$. In the previous experiment (Figure~\ref{fig:Experiment1}), a cell transfer of $10^6$ resulted in 8.5 cells/$10^5$ leukocytes measured in the lymph nodes on day 0. Doubling this, we set the T cell dose at 17 cells/$10^5$ leukocytes for each injection and supply it using $f_N(t)$ for CFSE-labelled cells and $f_{\tilde{N}}(t)$ for the competing cohort of cells. To match the experimental design, we set $f_{\tilde{N}}(t)$ as follows:
 \begin{itemize}
     \item Group (i) - control group, no competing cohort, $f_{\tilde{N}}(t) = 0$.
     \item Group (ii) - competing cohort injected at the same time as labelled cells ($\tilde{\mu}=27$).
     \item Group (iii) - competing cohort injected 24 hours prior to  labelled cells ($\tilde{\mu}=3$).
 \end{itemize}
 
\begin{figure}
    \centering
    \includegraphics[width=\textwidth]{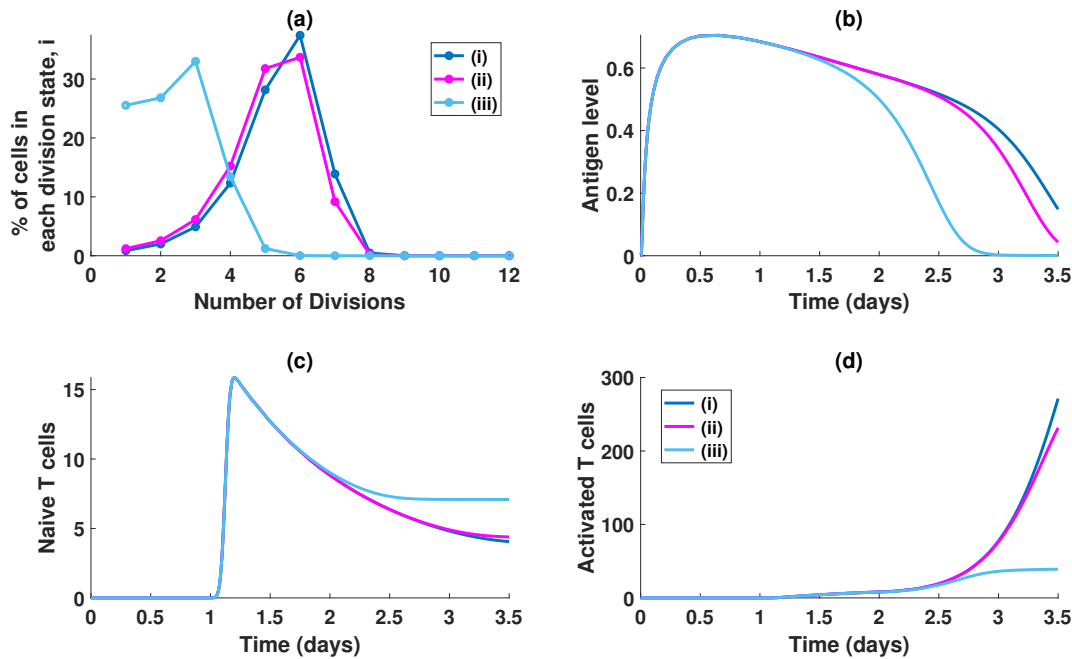}
    \caption{Simulations for experiment 2, Figure~\ref{fig:Experiment2}. In (a) there are no significant differences in division profiles if CFSE-labelled cells are injected on the same day as the competing cohort. A 24-hour delay in injection causes a reduction in divisions per cell. Results are in line with experimental observations, Figure~\ref{fig:Experiment2}d. The limited divisions experienced by group (iii) are explained in (b) where antigen depletion begins earlier due to the competing cohort injected at day 0. In (c) and (d), Total naive and activated T cell populations are very similar in groups (i) and (ii) but significantly different for group (iii).}
    \label{fig:Experiment 2-1}
\end{figure}

\begin{figure}
    \centering
    \includegraphics[width=\textwidth]{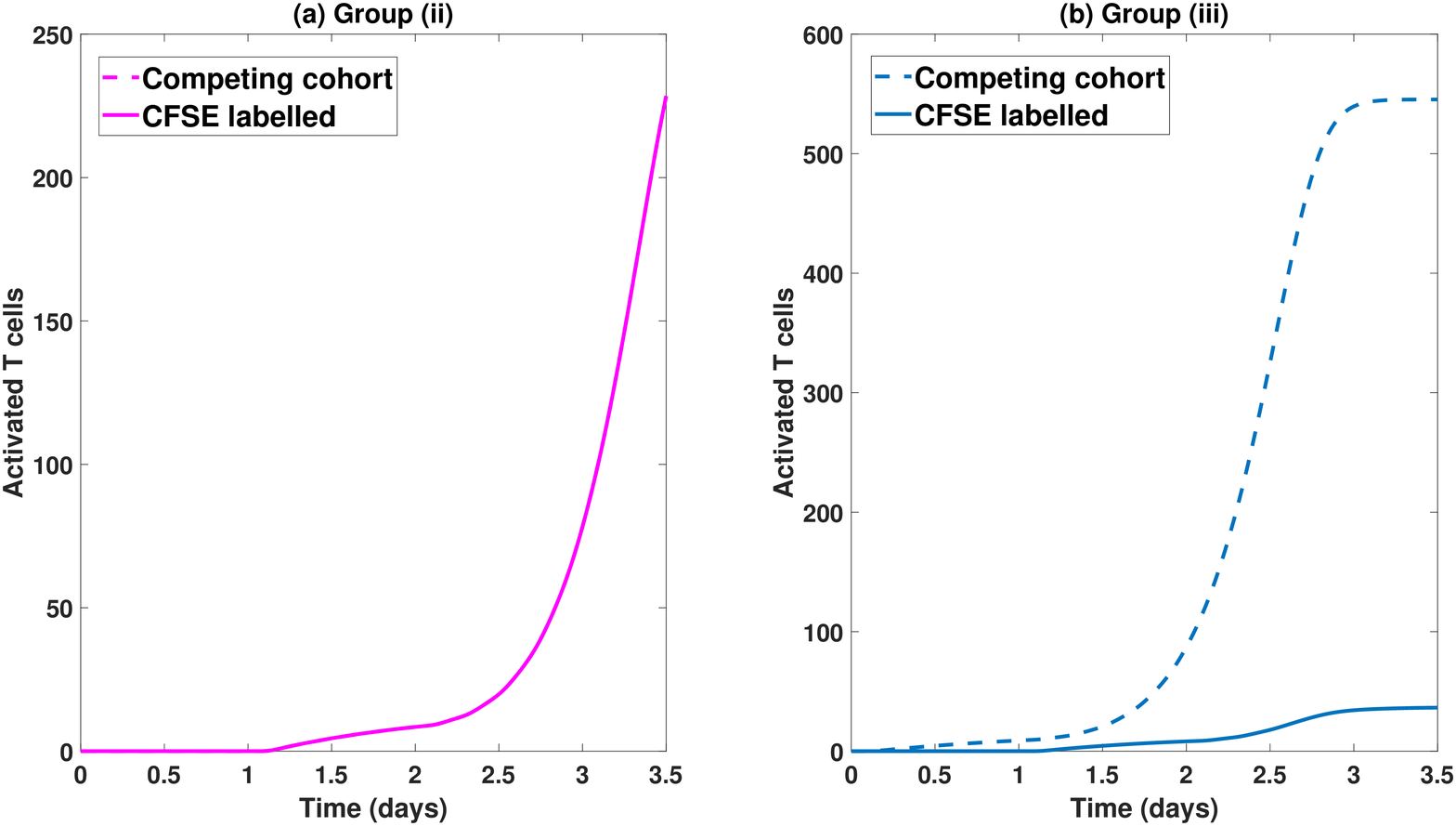}
    \caption{Total activated T cells from each cohort for (a) group (ii) and (b) group (iii). In (a) group (ii), when labelled cells are transferred at the same time as the competing cells, the populations contribute equally to the overall expansion (trajectories trace each other). In (b) group (iii), with a 24-hour delay between the competing cohort and labelled cells, the labelled cells only contribute minimally to the overall expansion.}
    \label{fig:Experiment 2-2}
\end{figure}

We now compare experimental results in Figure~\ref{fig:Experiment2} with model output in Figure~\ref{fig:Experiment 2-1}, tracking the CFSE-labelled cells. Our results in Figure~\ref{fig:Experiment 2-1}a demonstrate that if the competing cohort is released on the same day as the labelled cells (group (ii)), there is no change in the division profile compared with control group (i). However, if the competing cohort is injected one day prior to the labelled cells, the model estimates a significant drop in division profiles in line with experimental results, Figure~\ref{fig:Experiment2}c. In Figure~\ref{fig:Experiment 2-1}b, we see how antigen downregulation by T cells affects each population. No significant difference in loss of antigen is apparent between groups (i) and (ii) until approximately day 3. For group (iii), downregulation of antigen begins earlier, limiting the division of labelled cells as the competing cohort is injected at an earlier time. The similarity in antigen availability for groups (i) and (ii) means that naive and activated T cells experience almost identical dynamics in the timespan of the experiment as in Figures~\ref{fig:Experiment 2-1}c and \ref{fig:Experiment 2-1}d respectively. On the other hand, the delay of 24 hours in group (iii) causes the naive population to stabilise much earlier than when we have no delay (group (ii)) or no competing cohort (group (i)). As a result the total activated population in group (iii) reaches peak numbers around day 3 while the other groups continue to grow. Using equation (\ref{eq:recruitment}), with $N(0) = 17$, there is a minimal effect on cell recruitment into division between (i) and (ii) from $76\%$ to $74\%$. The delayed cohort, group (iii), only experiences $58\%$ recruitment. The estimated percentages are a good approximation to the experimental results in Figure~\ref{fig:Experiment2}b. We can compare the contribution of total activated T cells in the immune response from each cohort (competing cells versus labelled cells) in Figure~\ref{fig:Experiment 2-2}. From Figure~\ref{fig:Experiment 2-2}a, when both cohorts are injected at the same time, they contribute equally to the overall expansion (trajectories trace each other). With a 24-hour delay between transfer of competing cohort and labelled cells, there is a significant drop in the contribution of labelled cells. These results suggest that for a given number of initial precursors, the observed reduction in T cell divisions and recruitment can be explained by restrictions in antigen availability as growing cohorts of T cells downregulate antigen ensuring that T cells arriving late to an ongoing response will experience only limited activation and proliferation.\\


\begin{figure}
     \centering
     \includegraphics[width=\textwidth]{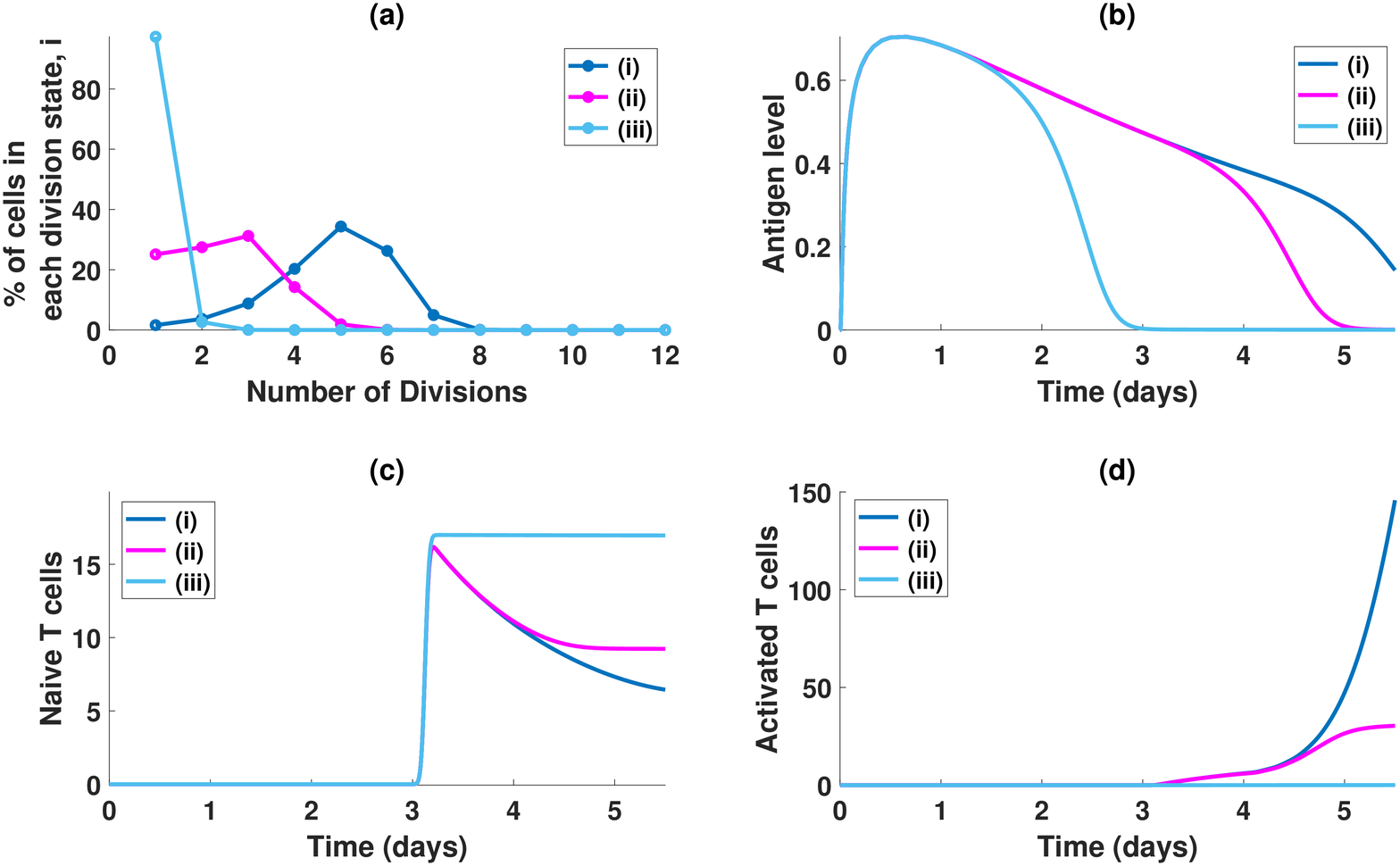}
     \caption{Simulations for Experiment 3, Figure~\ref{fig:Experiment3}. In (a), there is a significant change in number of divisions per recruited cell as the delay between competing cohort of T cells and labelled T cells is increased from control group (i) to a 24-hour delay in group (iii) and a 72-hour delay in group (iii). These division profiles are comparable with those in Figure~\ref{fig:Experiment3}d. In (b), the antigen level decreases more rapidly as the delay is increased, controlling T cell divisions. In (c) and (d), cell recruitment is drastically reduced in group (iii) with a 72-hour delay.}
     \label{fig:Experiment3-1}
 \end{figure}

\begin{figure}
     \centering
     \includegraphics[width=\textwidth]{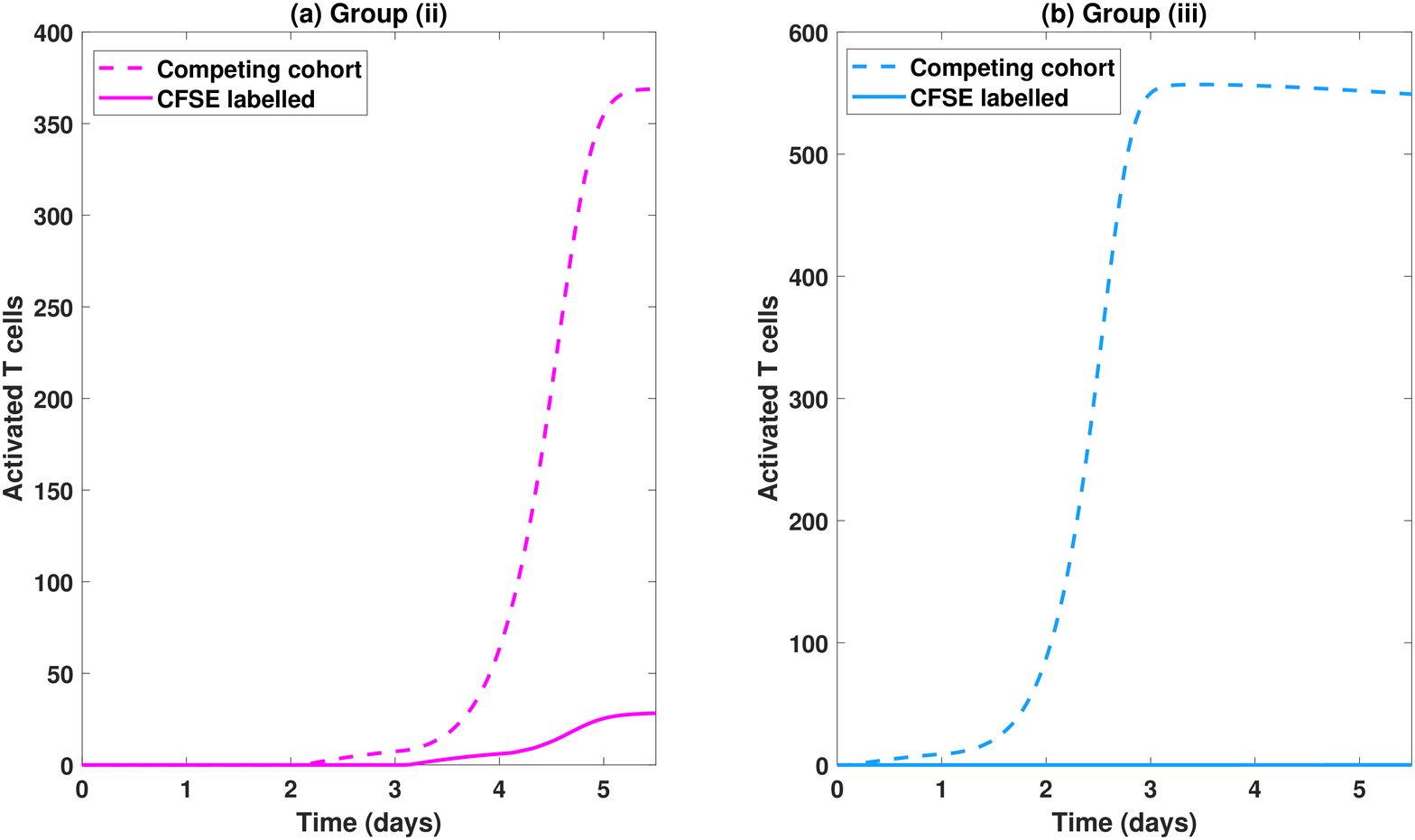}
     \caption{Activated T cells from each cohort for (a) group (ii) and (b) group (iii). In (a), labelled cells are transferred 24 hours following transfer of competing cohort. Labelled cells contribute minimally to the overall clonal expansion. With a 3-day delay in (b), the contribution of labelled cells is insignificant.}
     \label{fig:Experiment3-2}
 \end{figure}

\begin{figure}
    \centering
    \includegraphics[scale=0.5]{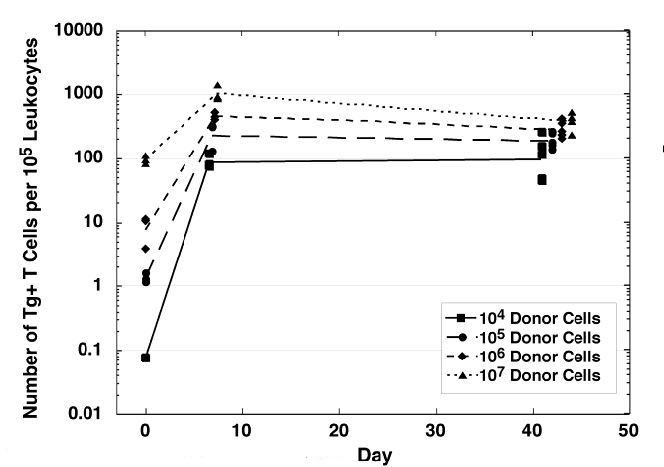}
    \caption{Experiment 1 results from Smith {\em et al.} \citep{Smith2020.09.09.290627}. The graph shows the number of Tg+ T cells per $10^5$ leukocytes measured in the draining lymph nodes (LN). Experimental procedure: Transgenic LN T cells are injected into non-transgenic Ly5.1 congenic recipient mice in doses of $10^4$, $10^5$, $10^6$ and $10^7$. Four days after transfer, 3 mice from each group are sacrificed to establish initial number of precursors in lymph nodes and spleen. Three days later, on day 0 of experiment, the remaining 8 mice from each group (32 in total) are immunised subcutaneously with $10\mu$g of moth cytochrome c (MCC$_{87-103}$). Draining LNs, distal LNs, spleens and peripheral blood is examined on days 7 and 42 after immunisation.}
    \label{fig:Experiment1}
\end{figure}

\begin{figure}
     \centering
     \includegraphics[scale=0.4]{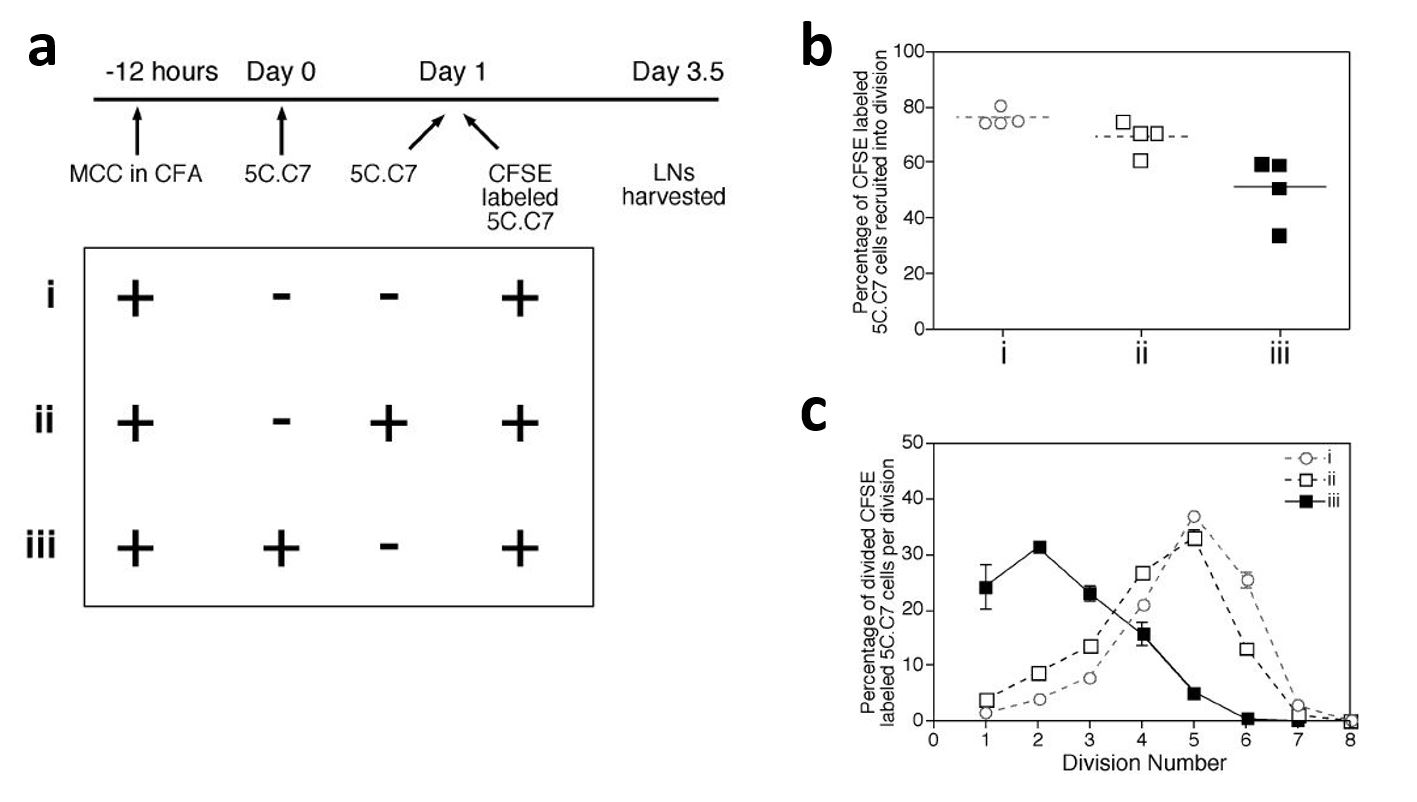}
     \caption{Experiment 2 plan and results from Spencer {\em et al.} \citep{spencer2020antigen}. \textbf{(a)} Experimental plan. Three groups of 4 recipient mice were immunised with 10$\mu$g MCC peptide in CFA and 1.5 days later were adoptively transferred with CFSE-labelled 5C.C7 TCR transgenic lymph node cells containing $2.5 \times 10^6$ CD4+ T cells. In the control group (i), no other cells were administered. Group (ii) received an equal number of 5C.C7 cells at the same time as the CFSE-labelled cohort, whereas group (iii) received an equal number of 5C.C7 cells 1 day before the CFSE-labelled cohort. Draining lymph nodes were harvested 2.5 days after the CFSE-labelled adoptive transfer. \textbf{(b)} Percentage of CFSE-labelled CD4+ TCR transgenic T cells recruited into cell division. Each point represents an individual mouse and the group mean is indicated by the horizontal line. \textbf{(c)} Percentage of divided CFSE-labelled CD4+ TCR transgenic T cells in each division peak at day 3.5. The number of cells in each division peak was divided by the total number of divided cells (i.e,, excluding undivided cells). Each point represents the mean calculated from the 4 mice in each group and the error bars show the standard error measurement.}
     \label{fig:Experiment2}
 \end{figure}

\begin{figure}
     \centering
     \includegraphics[scale=0.4]{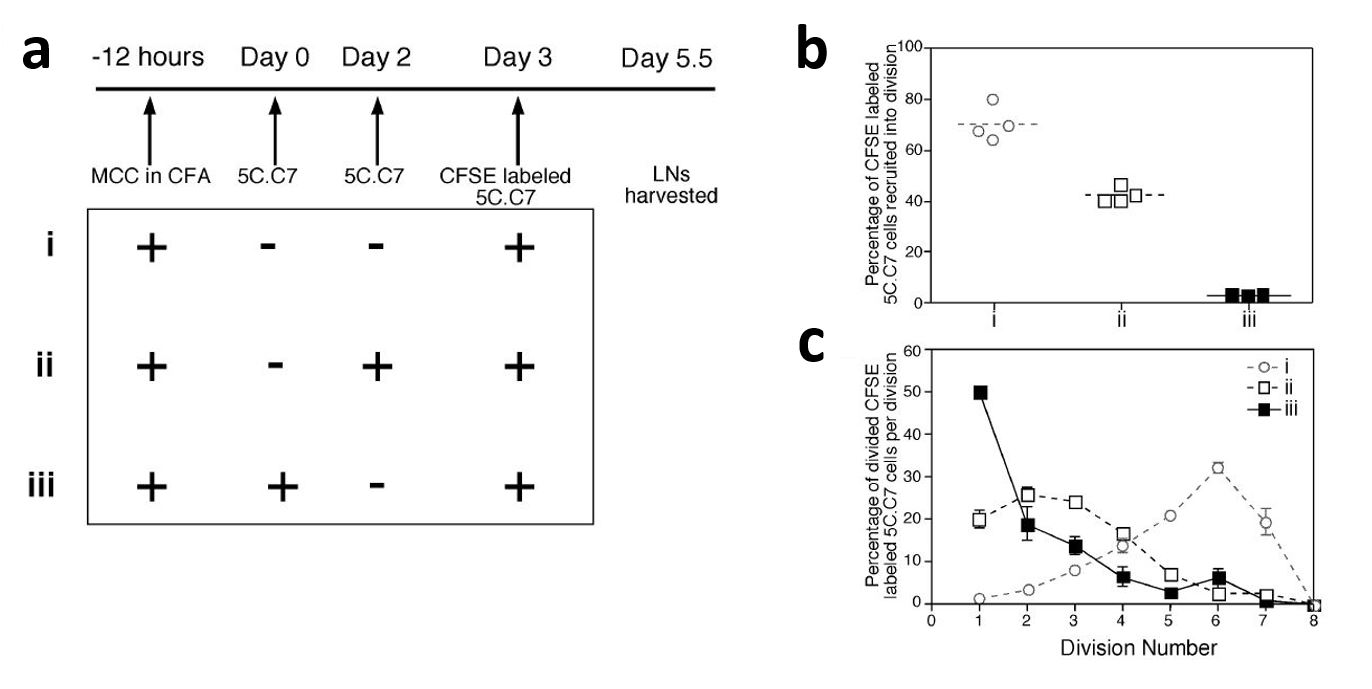}
     \caption{Experiment 3 plan and results from Spencer {\em et al.} \citep{spencer2020antigen}. \textbf{(a)} Experimental plan. Three groups of 4 recipient mice were immunised with 10$\mu$g MCC peptide in CFA and 3.5 days later were adoptively transferred with CFSE-labelled 5C.C7 TCR transgenic lymph node cells containing $2.5 \times 10^6$ CD4+ T cells. In the control group (i), no other cells were administered. Groups (ii) and (iii) received an equal number of 5C.C7 cells 1 or 3 days before the CFSE-labelled cohort. Draining lymph nodes were harvested 2.5 days after the CFSE-labelled adoptive transfer. \textbf{(b)} Percentage of CFSE-labelled TCR transgenic CD4+ T cells recruited into cell division. Each point represents an individual mouse and the group mean is indicated by the horizontal line. (d) Percentage of divided CFSE-labelled TCR transgenic CD4+ T cells in each division peak at day 5.5. Each point represents the mean calculated from the 4 mice in each group and the error bars show the standard error measurement.}
     \label{fig:Experiment3}
 \end{figure}

We now simulate the experimental set up in Figure~\ref{fig:Experiment3}, which tests the response of T cells to a longer delay between the labelled cells and competing cohort. Following the experimental design in Figure~\ref{fig:Experiment3}a, the CFSE-labelled cells are injected on day 3. We set $\mu = 75$ in $f_N(t)$ so that peak of T cell supply occurs at $t=75$ hours. We also need to specify the supply of the competing cohort. From experimental design, they are injected at $t=48$ hours for group (ii) and $t=0$ for group (iii). We set $f_{\tilde{N}}(t)$ as follows:
\begin{itemize}
     \item Group (i) - control group, no competing cohort, $f_{\tilde{N}}(t) = 0$.
     \item Group (ii) - competing cohort injected at t=48 ($\tilde{\mu}=51$).
     \item Group (iii) - competing cohort injected at t=0 ($\tilde{\mu}=3$).
 \end{itemize}
Our results indicate a significant drop in division profiles for the labelled cells as the delay in transfer of CFSE cells is increased from no delay (no competing cohort) in group (i), to a 1-day delay (competing cohort injected 1 day earlier) in group (ii) and finally a 3-day delay (competing cohort injected 3 days earlier) in group (iii), Figure~\ref{fig:Experiment3-1}a. Figure~\ref{fig:Experiment3-1}b shows that antigen downregulation by T cells is very rapid in group (iii) due to the the growing competing cohort injected at t=0. The effect of rapid downregulation of antigen on the labelled T cell populations can be seen in Figure~\ref{fig:Experiment3-1}(c-d), where hardly any naive T cells get activated in group (iii). Our model calculates the number of cells recruited into division as $62\%$, $46\%$ and $0.21\%$ for groups (i), (ii), and (iii) respectively which is inline with experimental outcome in Figure~\ref{fig:Experiment3}b. Comparing the total number of activated cells in each cohort, Figure~\ref{fig:Experiment3-2} indicates that increasing the delay in arrival of cells to an ongoing response significantly reduces the contribution of the CFSE-labelled cells.\\

Overall, the simulations from our model are a close approximation to experimental results showing that limitations in antigen supply, along with a reduction in T cell proliferation rate can explain the controlled behaviour of CD4+ T cell dynamics during an immune response. Specifically the model recreates the findings where clonal expansion of T cells reduce as precursor numbers increase as in Experiment 1 and T cell recruitment in to division and subsequent cell burst size drastically reduce as T cells are delayed into an ongoing immune response as in Experiments 2 and 3. We will now present a general discussion as to the significance of the results. 

\section{Discussion}\label{Discussion}
The experiments by Smith {\em et al.} and Spencer {\em et al.} demonstrated that the immune system exerts a tight control on CD4+ T cell clonal expansion by restricting recruitment into division and the number of divisions experienced per parent cell. To explain this behaviour, we used a mathematical model describing a negative feedback between antigen stimulation of T cells and T cell downregulation of antigen. In this way, we simulated an environment where T cells compete for limited antigen. Fitting our model to Experiment 1 demonstrated that antigen availability can explain the reduction in T cell expansion as initial precursors were increased by factors of ten. This suggests that the clonal expansion during an infection may be determined and limited by the magnitude of infection and antigen downregulation. Interestingly, we found that when antigen is available for longer, as with the smallest precursor population, a reduction in T cell proliferation rate may be responsible in further controlling the clonal expansion. This extra control would insure that cell division is limited at high antigen doses. Limitations in clonal expansion when antigen is increased has been observed experimentally. An {\em in vitro} study by Jelley-Gibbs {\em et al.} \citep{jelley2005repeated} showed that when antigen presentation is limited to 2 days, maximum CD4+ effector expansion is achieved while repeated antigen stimulation impaired T cell generation. Quiel {\em et al.} \citep{quiel2011antigen} measured the factor of expansion (ratio of the number of antigen-specific T cells at day 7 to the number of cells before immunisation) in two different sizes of antigen-specific precursors, 300 and 30\ 000 cells, when the cells were exposed to two levels of antigen, 100 $\mu$g and 1 mg. They found that increasing antigen dose resulted in a modest increase in the factor of expansion in each population size (\citep{quiel2011antigen}, Figure 5A). The reduction in the proliferation rate as T cells continue to be stimulated by antigen could also arise from T cell anergy~\citep{schwartz2003t}. This is a process in which T cells become unresponsive to the proliferative effect of IL-2 following antigen exposure. Yamamoto et al. demonstrated in vivo, that Th1, CD4+ T cells become anergic depending on intensity of antigen stimulation and duration of exposure~\citep{yamamoto2007induction}. Therefore, if a proportion of T cells become anergic after a period of time, then the overall proliferation rate would reduce.\\

Our model estimated a $6\%$ drop in T cell recruitment into division for each factor 10 increase in precursors, so that when a competing cohort of cells was released at the same time as the labelled cells, effectively doubling the initial T cell dose (Experiment 2, group (ii)), there was no significant change in
recruitment into division or burst size per cell. This can be seen by comparing the division profiles for control group (i) with no competing cohort against group (ii) as in Figure~\ref{fig:Experiment 2-1}a. Our result captures the experimental outcome showing that there is no significant alteration in recruitment of the CFSE-labelled cells when there is no delay in response or if the response is in the same order of magnitude. Tracking activation profiles, we found that the two populations contributed equally to the response.\\ 

Increasing the delay between the competing cohort and the labelled cells reduced both recruitment and burst size in line with experimental results. The activation profiles from these simulations revealed that cells arriving late to an ongoing response were excluded from recruitment and division, and this exclusion increased for increasing delay (from 24 hours to 72 hours). Our model attributed this drastic change to antigen downregulation by cells arriving earlier to the response. Mayerova et al. obtained similar results, showing a decrease in expansion of a delayed cohort for CD8+ T cells, \citep{mayerova2006conditioning}. In their case, the authors attributed the reduction in clonal expansion to phenotypic changes in the APCs (Langerhans cells), which were proposed to compromise the stimulation of the naïve T cells; however, additional experiments by Spencer et al, showed that a competing cohort responding to the same APCs but an unrelated antigen specificity had no effect on the activation of labelled cells, regardless of delay \citep{spencer2020antigen}. This allowed the authors to exclude changes in APC phenotype and other antigen non-specific effects (such as T cell crowding around APCs) as causes of reduction in recruitment and division.

Studies have also suggested that cytokine signalling molecules such as interleukins can control the response of T cells to antigen presentation. One such cytokine is interleukin 2 (IL-2), a cytokine linked to T cell growth, differentiation and suppression of the immune response \citep{bachmann2007interleukin}. Experiments by Villarino {\em et al.} suggested a negative feedback loop for the production of IL-2 by CD4+ T cells \citep{villarino2007helper}. They demonstrated {\em in vivo}, that when CD4+ T cells are introduced into an environment rich in IL-2, they suppress their production of this growth signal during an immune response. This effect was deemed to be transient and requires continual availability of IL-2. In another study by Blattman {\em et al.} the introduction of IL-2 during the expansion phase reduced peak numbers of CD4+ T cells by $80\%$, interestingly IL-2 introduction did not affect CD8+ T-cell response \citep{blattman2003therapeutic}. In contrast, IL-2 introduction during the contraction and memory phases resulted in increased proliferation and survival of T cells. IL-12 and interferon-gamma (IFN-$\gamma$) have also been shown to impact T cell response. Eriksson {\em et al.} showed that CD4+ cell proliferation in IL-12 deficient mice was poor \citep{eriksson2001dual}. On the other hand, mice that were deficient in IFN-$\gamma$ developed autoimmune disease. The role of IFN-$\gamma$ as a cytokine responsible for the downregulation of T-cell response has been reported in other studies \citep{badovinac2000regulation, hosking2016tcr, eriksson2001lethal, bachmann2002balancing}. Regulatory T cells (Tregs) have also been implicated in the suppression of T-cell response. Dowling {\em et al.} demonstrated that Tregs can mediate the magnitude of effector T cells by changing the number of overall T cell divisions \citep{dowling2018regulatory}. In an opinion article by Gasteiger and Kastenmuller, it was proposed that the availability of IL-2 is mediated by Tregs, restricting the CD4+ t cell response \citep{gasteiger2012foxp3+}.\\

The examples above demonstrate that the complexity of the interactions during a T-cell mediated immune response has led to a diverse range of experimental results. A variety of mathematical models also explain these dynamics by investigating the affects of cytokines and Tregs. As an example, Ganusov {\em et al.} produced a mathematical model that describes the in vitro dynamics of CD4+ T cells with respect to IL-2 concentration \citep{ganusov20072}. The model suggests that at low IL-2 concentrations fewer cells are recruited into division and the rate of cell death also increases. Overall, the authors summarised that IL-2 concentration affects cell death rate and not cell division rate. An overview of combined mathematical models and experiments can be found in \citep{morel2014modeling}. We leave this discussion by considering that it is likely that T cells exhibit many traits during an immune response and in modelling the dynamics we can attempt to find the extent of dependency of the clonal expansion on these characteristics. In our model, the clonal expansion of T cells can successfully be described by antigen downregulation, a conclusion that is supported by the experimental findings that competing responses to a different antigen have no effect on clonal expansion~\citep{spencer2020antigen}; however, the addition of cytokines and Tregs into the model could assist in understanding the impact of signalling in this scenario and is an opportunity for future investigation. 

\section{Conclusion}

Our results indicate that the dynamics of T cell expansion may heavily depend on antigen availability. From this perspective, T cells require continual communication with APCs during the expansion phase leading to an antigen competitive environment. Subsequently, competition for antigen leads to an appropriate immune response corresponding to the magnitude of infection. Further to this, our results also suggest that an extra layer of control is established by an intrinsic reduction in T cell proliferation rate, especially for a large antigen dose or prolonged antigen exposure. Testing our model against experimental data where the extent of T cell expansion is measured in relation to variations in antigen dose and duration of infection will greatly assist in calibrating this model to gain a deeper understanding into the mechanisms that drive a robust and tailored response to pathogens.\\

\begin{acknowledgements}
The authors gratefully acknowledge support for this work through Australian Government Research Training Program Scholarship (PP) and the Australian Research Council Discovery Project DP180101512 (PSK).
\end{acknowledgements}

\clearpage
\bibliographystyle{spbasic}      
\bibliography{mybibfile}   

\begin{thebibliography}{40}
\providecommand{\natexlab}[1]{#1}
\providecommand{\url}[1]{{#1}}
\providecommand{\urlprefix}{URL }
\expandafter\ifx\csname urlstyle\endcsname\relax
  \providecommand{\doi}[1]{DOI~\discretionary{}{}{}#1}\else
  \providecommand{\doi}{DOI~\discretionary{}{}{}\begingroup
  \urlstyle{rm}\Url}\fi
\providecommand{\eprint}[2][]{\url{#2}}

\bibitem[{Allan et~al(2006)Allan, Waithman, Bedoui, Jones, Villadangos, Zhan,
  Lew, Shortman, Heath, and Carbone}]{allan2006migratory}
Allan RS, Waithman J, Bedoui S, Jones CM, Villadangos JA, Zhan Y, Lew AM,
  Shortman K, Heath WR, Carbone FR (2006) Migratory dendritic cells transfer
  antigen to a lymph node-resident dendritic cell population for efficient ctl
  priming. Immunity 25(1):153--162

\bibitem[{Bachmann and Kopf(2002)}]{bachmann2002balancing}
Bachmann MF, Kopf M (2002) Balancing protective immunity and immunopathology.
  Curr Opin Immunol 14(4):413--419

\bibitem[{Bachmann and Oxenius(2007)}]{bachmann2007interleukin}
Bachmann MF, Oxenius A (2007) Interleukin 2: from immunostimulation to
  immunoregulation and back again. EMBO Rep 8(12):1142--1148

\bibitem[{Badovinac et~al(2000)Badovinac, Tvinnereim, and
  Harty}]{badovinac2000regulation}
Badovinac VP, Tvinnereim AR, Harty JT (2000) Regulation of antigen-specific
  cd8+ t cell homeostasis by perforin and interferon-$\gamma$. Science
  290(5495):1354--1357

\bibitem[{Beverley and Maini(2000)}]{beverley2000differences}
Beverley PC, Maini MK (2000) Differences in the regulation of cd4 and cd8
  t--cell clones during immune responses. Philos Trans R Soc, B
  355(1395):401--406

\bibitem[{Blattman et~al(2003)Blattman, Grayson, Wherry, Kaech, Smith, and
  Ahmed}]{blattman2003therapeutic}
Blattman JN, Grayson JM, Wherry EJ, Kaech SM, Smith KA, Ahmed R (2003)
  Therapeutic use of il-2 to enhance antiviral t-cell responses in vivo. Nat
  Med 9(5):540--547

\bibitem[{Borghans et~al(1999)Borghans, Taams, Wauben, and
  De~Boer}]{borghans1999competition}
Borghans JA, Taams LS, Wauben MH, De~Boer RJ (1999) Competition for antigenic
  sites during t cell proliferation: a mathematical interpretation of in vitro
  data. Proc Natl Acad Sci U S A 96(19):10,782--10,787

\bibitem[{den Braber et~al(2012)den Braber, Mugwagwa, Vrisekoop, Westera,
  M{\"o}gling, de~Boer, Willems, Schrijver, Spierenburg, Gaiser
  et~al}]{den2012maintenance}
den Braber I, Mugwagwa T, Vrisekoop N, Westera L, M{\"o}gling R, de~Boer AB,
  Willems N, Schrijver EH, Spierenburg G, Gaiser K, et~al (2012) Maintenance of
  peripheral naive t cells is sustained by thymus output in mice but not
  humans. Immunity 36(2):288--297

\bibitem[{De~Boer and Perelson(2013)}]{de2013antigen}
De~Boer RJ, Perelson AS (2013) Antigen-stimulated cd4 t cell expansion can be
  limited by their grazing of peptide--mhc complexes. J Immunol
  190(11):5454--5458

\bibitem[{Dhainaut and Moser(2014)}]{dhainaut2014regulation}
Dhainaut M, Moser M (2014) Regulation of immune reactivity by intercellular
  transfer. front immunol 5: 112

\bibitem[{Diao et~al(2006)Diao, Winter, Cantin, Chen, Xu, Kelvin, Phillips, and
  Cattral}]{diao2006situ}
Diao J, Winter E, Cantin C, Chen W, Xu L, Kelvin D, Phillips J, Cattral MS
  (2006) In situ replication of immediate dendritic cell (dc) precursors
  contributes to conventional dc homeostasis in lymphoid tissue. J Immunol
  176(12):7196--7206

\bibitem[{Dowling et~al(2018)Dowling, Kan, Heinzel, Marchingo, Hodgkin, and
  Hawkins}]{dowling2018regulatory}
Dowling MR, Kan A, Heinzel S, Marchingo JM, Hodgkin PD, Hawkins ED (2018)
  Regulatory t cells suppress effector t cell proliferation by limiting
  division destiny. Front Immunol 9:2461

\bibitem[{Eriksson et~al(2001{\natexlab{a}})Eriksson, Kurrer, Bingisser,
  Eugster, Saremaslani, Follath, Marsch, and Widmer}]{eriksson2001lethal}
Eriksson U, Kurrer M, Bingisser R, Eugster H, Saremaslani P, Follath F, Marsch
  S, Widmer U (2001{\natexlab{a}}) Lethal autoimmune myocarditis in
  interferon-$\gamma$ receptor--deficient mice: Enhanced disease severity by
  impaired inducible nitric oxide synthase induction. Circulation 103(1):18--21

\bibitem[{Eriksson et~al(2001{\natexlab{b}})Eriksson, Kurrer, Sebald,
  Brombacher, and Kopf}]{eriksson2001dual}
Eriksson U, Kurrer MO, Sebald W, Brombacher F, Kopf M (2001{\natexlab{b}}) Dual
  role of the il-12/ifn-$\gamma$ axis in the development of autoimmune
  myocarditis: induction by il-12 and protection by ifn-$\gamma$. J Immunol
  167(9):5464--5469

\bibitem[{Furuta et~al(2012)Furuta, Ishido, and Roche}]{furuta2012encounter}
Furuta K, Ishido S, Roche PA (2012) Encounter with antigen-specific primed cd4
  t cells promotes mhc class ii degradation in dendritic cells. Proc Natl Acad
  Sci U S A 109(47):19,380--19,385

\bibitem[{Ganusov et~al(2007)Ganusov, Milutinovi{\'c}, and
  De~Boer}]{ganusov20072}
Ganusov VV, Milutinovi{\'c} D, De~Boer RJ (2007) Il-2 regulates expansion of
  cd4+ t cell populations by affecting cell death: insights from modeling cfse
  data. J Immunol 179(2):950--957

\bibitem[{Gasteiger and Kastenmuller(2012)}]{gasteiger2012foxp3+}
Gasteiger G, Kastenmuller W (2012) Foxp3+ regulatory t-cells and il-2: the
  moirai of t-cell fates? Front Immunol 3:179

\bibitem[{Homann et~al(2001)Homann, Teyton, and
  Oldstone}]{homann2001differential}
Homann D, Teyton L, Oldstone MB (2001) Differential regulation of antiviral
  t-cell immunity results in stable cd8+ but declining cd4+ t-cell memory. Nat
  Med (N Y, NY, U S) 7(8):913

\bibitem[{Hosking et~al(2016)Hosking, Flynn, and Whitton}]{hosking2016tcr}
Hosking MP, Flynn CT, Whitton JL (2016) Tcr independent suppression of cd8+ t
  cell cytokine production mediated by ifn$\gamma$ in vivo. Virology 498:69--81

\bibitem[{Jelley-Gibbs et~al(2005)Jelley-Gibbs, Dibble, Filipson, Haynes, Kemp,
  and Swain}]{jelley2005repeated}
Jelley-Gibbs DM, Dibble JP, Filipson S, Haynes L, Kemp RA, Swain SL (2005)
  Repeated stimulation of cd4 effector t cells can limit their protective
  function. J Exp Med 201(7):1101--1112

\bibitem[{Kamath et~al(2002)Kamath, Henri, Battye, Tough, and
  Shortman}]{kamath2002developmental}
Kamath AT, Henri S, Battye F, Tough DF, Shortman K (2002) Developmental
  kinetics and lifespan of dendritic cells in mouse lymphoid organs. Blood
  100(5):1734--1741

\bibitem[{Kedl et~al(2000)Kedl, Rees, Hildeman, Schaefer, Mitchell, Kappler,
  and Marrack}]{kedl2000t}
Kedl RM, Rees WA, Hildeman DA, Schaefer B, Mitchell T, Kappler J, Marrack P
  (2000) T cells compete for access to antigen-bearing antigen-presenting
  cells. J Exp Med 192(8):1105--1114

\bibitem[{Lambrecht et~al(2000)Lambrecht, Pauwels, and
  Groth}]{lambrecht2000induction}
Lambrecht BN, Pauwels RA, Groth BFdS (2000) Induction of rapid t cell
  activation, division, and recirculation by intratracheal injection of
  dendritic cells in a tcr transgenic model. J Immunol 164(6):2937--2946

\bibitem[{Liou et~al(2012)Liou, Myers, Barkauskas, and
  Huang}]{liou2012intravital}
Liou HR, Myers JT, Barkauskas DS, Huang AY (2012) Intravital imaging of the
  mouse popliteal lymph node. J Visualized Exp (60):e3720

\bibitem[{Mayer et~al(2019)Mayer, Zhang, Perelson, and
  Wingreen}]{mayer2019regulation}
Mayer A, Zhang Y, Perelson AS, Wingreen NS (2019) Regulation of t cell
  expansion by antigen presentation dynamics. Proc Natl Acad Sci U S A
  116(13):5914--5919

\bibitem[{Mayerova et~al(2006)Mayerova, Wang, Bursch, and
  Hogquist}]{mayerova2006conditioning}
Mayerova D, Wang L, Bursch LS, Hogquist KA (2006) Conditioning of langerhans
  cells induced by a primary cd8 t cell response to self-antigen in vivo. J
  Immunol 176(8):4658--4665

\bibitem[{Morel et~al(2014)Morel, Faeder, Hawse, and
  Miskov-Zivanov}]{morel2014modeling}
Morel PA, Faeder JR, Hawse WF, Miskov-Zivanov N (2014) Modeling the t cell
  immune response: a fascinating challenge. J Pharmacokinet Pharmacodyn
  41(5):401--413

\bibitem[{Obst(2015)}]{obst2015timing}
Obst R (2015) The timing of t cell priming and cycling. Front Immunol 6:563

\bibitem[{Pappalardo et~al(2016)Pappalardo, Fichera, Paparone, Lombardo,
  Pennisi, Russo, Leotta, Pappalardo, Pedretti, De~Fiore
  et~al}]{pappalardo2016computational}
Pappalardo F, Fichera E, Paparone N, Lombardo A, Pennisi M, Russo G, Leotta M,
  Pappalardo F, Pedretti A, De~Fiore F, et~al (2016) A computational model to
  predict the immune system activation by citrus-derived vaccine adjuvants.
  Bioinformatics 32(17):2672--2680

\bibitem[{Pennisi et~al(2019)Pennisi, Russo, Sgroi, Bonaccorso, Palumbo,
  Fichera, Mitra, Walker, Cardona, Amat et~al}]{pennisi2019predicting}
Pennisi M, Russo G, Sgroi G, Bonaccorso A, Palumbo GAP, Fichera E, Mitra DK,
  Walker KB, Cardona PJ, Amat M, et~al (2019) Predicting the artificial
  immunity induced by ruti{\textregistered} vaccine against tuberculosis using
  universal immune system simulator (uiss). BMC bioinform 20(6):1--10

\bibitem[{Quiel et~al(2011)Quiel, Caucheteux, Laurence, Singh, Bocharov,
  Ben-Sasson, Grossman, and Paul}]{quiel2011antigen}
Quiel J, Caucheteux S, Laurence A, Singh NJ, Bocharov G, Ben-Sasson SZ,
  Grossman Z, Paul WE (2011) Antigen-stimulated cd4 t-cell expansion is
  inversely and log-linearly related to precursor number. Proc Natl Acad Sci U
  S A 108(8):3312--3317

\bibitem[{Rabenstein et~al(2014)Rabenstein, Behrendt, Ellwart, Naumann, Horsch,
  Beckers, and Obst}]{rabenstein2014differential}
Rabenstein H, Behrendt AC, Ellwart JW, Naumann R, Horsch M, Beckers J, Obst R
  (2014) Differential kinetics of antigen dependency of cd4+ and cd8+ t cells.
  J Immunol 192(8):3507--3517

\bibitem[{Schwartz(2003)}]{schwartz2003t}
Schwartz RH (2003) T cell anergy. Annu Rev Immunol 21(1):305--334

\bibitem[{Smith and Fazekas~de St~Groth(2020)}]{Smith2020.09.09.290627}
Smith AL, Fazekas~de St~Groth B (2020) T cell competition profoundly reduces
  the effect of initial precursor frequency on the generation of cd4 t cell
  memory. bioRxiv \doi{10.1101/2020.09.09.290627},
  \urlprefix\url{https://www.biorxiv.org/content/early/2020/09/10/2020.09.09.290627},
  \eprint{https://www.biorxiv.org/content/early/2020/09/10/2020.09.09.290627.full.pdf}

\bibitem[{Spencer et~al(2020)Spencer, Smith, and
  de~St~Groth}]{spencer2020antigen}
Spencer AJ, Smith AL, de~St~Groth BF (2020) Antigen-specific competitive
  inhibition of cd4+ t cell recruitment into the primary immune response.
  bioRxiv

\bibitem[{Tomura et~al(2014)Tomura, Hata, Matsuoka, Shand, Nakanishi, Ikebuchi,
  Ueha, Tsutsui, Inaba, Matsushima et~al}]{tomura2014tracking}
Tomura M, Hata A, Matsuoka S, Shand FH, Nakanishi Y, Ikebuchi R, Ueha S,
  Tsutsui H, Inaba K, Matsushima K, et~al (2014) Tracking and quantification of
  dendritic cell migration and antigen trafficking between the skin and lymph
  nodes. Sci Rep 4:6030

\bibitem[{Van~de Velde and Murray(2016)}]{van2016proliferating}
Van~de Velde LA, Murray PJ (2016) Proliferating helper t cells require
  rictor/mtorc2 complex to integrate signals from limiting environmental amino
  acids. J Biol Chem 291(50):25,815--25,822

\bibitem[{Villarino et~al(2007)Villarino, Tato, Stumhofer, Yao, Cui,
  Hennighausen, O'Shea, and Hunter}]{villarino2007helper}
Villarino AV, Tato CM, Stumhofer JS, Yao Z, Cui YK, Hennighausen L, O'Shea JJ,
  Hunter CA (2007) Helper t cell il-2 production is limited by negative
  feedback and stat-dependent cytokine signals. J Exp Med 204(1):65--71

\bibitem[{Yamamoto et~al(2007)Yamamoto, Hattori, and
  Yoshida}]{yamamoto2007induction}
Yamamoto T, Hattori M, Yoshida T (2007) Induction of t-cell activation or
  anergy determined by the combination of intensity and duration of t-cell
  receptor stimulation, and sequential induction in an individual cell.
  Immunology 121(3):383--391

\bibitem[{Yarke et~al(2008)Yarke, Dalheimer, Zhang, Catron, Jenkins, and
  Mueller}]{yarke2008proliferating}
Yarke CA, Dalheimer SL, Zhang N, Catron DM, Jenkins MK, Mueller DL (2008)
  Proliferating cd4+ t cells undergo immediate growth arrest upon cessation of
  tcr signaling in vivo. J Immunol 180(1):156--162

\end{thebibliography}

\end{document}